# Решение линейных дробно-производных обыкновенных дифференциальных уравнений с постоянными матричными коэффициентами

**Алиев Фикрет А., Алиев Н.А., Сафарова Н.А., Касимова К.Г., Велиева Н.И.**

**Abstract.** Рассматривается задача Коши для дробно-производных линейных систем обыкновенных дифференциальных уравнений с постоянными коэффициентами, где впервые приводятся аналитические выражения через матричную экспоненту соответствующего ее решения. На основе полученных результатов приводятся условия, обеспечивающие асимптотическую устойчивость исходной системы. Результаты иллюстрируются на числовом примере, где показываются, что когда порядок дробной производной стремится к единице, то решение стремится к соответствующей экспоненциальной функции.

**Ключевые слова** Дробно-производные, постоянные матричные коэффициенты, обыкновенные дифференциальные уравнения, функция Миттаг-Лефлера, экспоненциальная функция.

## 1. Введение.

В последнее время много внимания уделяется исследованию дифференциальных уравнений с дробными производными в смысле Римана Лиувилля [1-4], в основном вопросам существования, единственности ее решения и т. п. Однако, к настоящему времени самый простой случай-решение системы линейных обыкновенных дифференциальных уравнений с постоянными коэффициентами дробной производной в явном виде не представлено и поэтому имеет смысл рассмотреть задачу Коши для этой системы. Также интересно представить решение задачи Коши от функции Миттаг-Лефлера к экспоненциальной функции, которое позволит вести анализ асимптотической устойчивости. Отметим, что в работе [5] эти задачи решены для одномерного случая и задача Коши для операторных обыкновенных дифференциальных уравнений. Также в [6] рассмотрена задача аналитического конструирования оптимального регулятора, когда движение объекта описывается линейной системой обыкновенных дифференциальных уравнений дробного порядка 1/3.

В данной работе рассматривается задача Коши для линейных систем обыкновенных дифференциальных уравнений дробного порядка с постоянными матричными коэффициентами и получена аналитическая формула для ее решения. Показано, что когда дробный порядок равняется единице, полученные в данной работе результаты совпадают с классическими результатами [7.8]. В отличие от [9] такой подход позволит найти аналитическое дискретное соотношение между любыми соседними точками. Результаты иллюстрируются числовыми примерами, где показывается что, когда порядок стремится к единице, то решение приближается к экспоненциальной функции.

## 2. Постановка задачи.

Рассмотрим следующую задачу Коши:

$$D^\alpha x(t) = Ax(t), \quad t > t_0, \quad x(t_0) = x_0, \qquad (1)$$

где $A$ $n \times n-$ мерная известная постоянная матрица, $x(t)$ $n-$мерный неизвестный вектор-столбец, $x_0$ -заданный $n-$мерный постоянный вектор, известное вещественное число $\alpha \in (0,1)$. Требуется найти вектор-столбец $x(t)$, который удовлетворяет задаче Коши (1).

Как известно [10], существует невырожденная матрица $T$, которая приводит матрицу $A$ к следующему каноническому виду

$$T^{-1}AT = \Lambda, \qquad (2)$$

где, в зависимости от характера собственных значений матрица $A$, матрица $\Lambda$ может быть или в диагональном, или в Жордановом или же в «комплексном» виде.

Для простоты рассмотрим случай, когда собственные значения матрицы $A$ вещественные и различные. Далее покажем, что оставшиеся два случая сводятся к приведенному с помощью метода продолжения коэффициента с малым параметром [11].

**3. Случай простых вещественных различных собственных значений матрица A.** Когда собственные значения матрицы A различно – вещественные, тогда матрица $\Lambda$ будет в следующем виде:

$$\Lambda = diag[\lambda_1, \lambda_2, ..., \lambda_n] \qquad (3)$$

где $\lambda_i (i = \overline{1.n})$ собственные значения матрицы $A$.

Учитывая (2) в (1) и обозначая $\mathrm{T}^{-1}x = y$ можно уравнение (1) переписать в следующем виде.

$$D^\alpha y = \Lambda y, \, y(t_0) = \mathrm{T}^{-1}x_0 \equiv \beta_0. \qquad (4)$$

Обозначая $y(t) = [y_1(t), y_2(t), ... y_n(t)]'$, где штрих означает операцию транспонирования, то система (4) переходит к виду

$$D^\alpha y_i(t) = \lambda_i y_i(t), \quad i = \overline{1.n}, \qquad (5)$$

где в [5] приведено ее решение в следующем компактном виде

$$y_i(t) = \beta_0 \frac{\sum_{s=0}^{2q} \lambda_i^{\frac{s}{2p+1}} \dfrac{d}{dt} \int_{t_0}^{t} \dfrac{(t-\tau)^{\frac{s-2q}{2q+1}}}{\dfrac{s-2q}{2q+1}!} e^{\tau \lambda_i^{\frac{2q+1}{2p+1}}} d\tau}{\left. \sum_{s=0}^{2q} \lambda_i^{\frac{s}{2p+1}} \dfrac{d}{dt} \int_{t_0}^{t} \dfrac{(t-\tau)^{\frac{s-2q}{2q+1}}}{\dfrac{s-2q}{2q+1}!} e^{\tau \lambda_i^{\frac{s-2q}{2q+1}}} d\tau \right|_{t=t_0}} \qquad (6)$$

Здесь $\alpha \approx \dfrac{2p+1}{2q+1}$, $p$ и q натуральные числа, т.е. $p, q \in N$. Как видно из (6) под знаком интеграла s<2q, а это приводит к неограниченности подынтегральной функции. Теперь преобразуем знаменатель выражения (6) так, чтобы под знаком интеграла при выполнении дифференцирования получилось ограниченное выражение. Для этого сперва подынтегральное выражение преобразуем так, чтобы получилась положительная степень. $(t-\tau)$ и затем интегрируя его по частям, уже можно будет дифференцировать по $\dfrac{d}{dt}$, т.е.

$$\sum_{s=0}^{2q}\lambda_i^{\frac{s}{2p+1}}\frac{d}{dt}\int_{t_0}^{t}\frac{(t-\tau)^{\frac{s-2q}{2q+1}}}{\frac{s-2q}{2q+1}!}e^{\tau\lambda_i^{\frac{2q+1}{2p+1}}}\,d\tau\Bigg|_{t=t_0}=\lambda_i^{\frac{2q}{2p+1}}\cdot e^{t_0\cdot\lambda_i^{\frac{2q+1}{2p+1}}} \tag{7}$$

Таким образом, учитывая (7) в (6) получим:

$$y_i(t)=\lambda_i^{-\frac{2q}{2p+1}}e^{-t_0\lambda_i^{\frac{2q+1}{2p+1}}}\cdot y_i(t_0)\cdot\sum_{s=0}^{2q}\lambda_i^{\frac{s}{2p+1}}\frac{d}{dt}\int_{t_0}^{t}\frac{(t-\tau)^{\frac{s-2q}{2q+1}}}{\frac{s-2q}{2p+1}!}e^{\tau\lambda_i^{\frac{2q+1}{2p+1}}}d\tau=$$

$$=\left[\sum_{s=0}^{2q}\lambda_i^{\frac{s-2q}{2p+1}}\cdot\frac{d}{dt}\int_{t_0}^{t}\frac{(t-\tau)^{\frac{s-2q}{2q+1}}}{\frac{s-2q}{2q+1}!}e^{(\tau-t_0)\lambda_i^{\frac{2q+1}{2p+1}}}d\tau\right]y_i(t_0) \tag{8}$$

Используя (8) для y(t) задачи Коши (1) имеем

$$y(t)=\left[\sum_{s=0}^{2q}\Lambda^{\frac{s-2q}{2p+1}}\frac{d}{dt}\int_{t_0}^{t}\frac{(t-\tau)^{\frac{s-2q}{2q+1}}}{\frac{s-2q}{2q+1}!}e^{(\tau-t_0)\Lambda^{\frac{2q+1}{2p+1}}}d\tau\right]y(t_0). \tag{9}$$

Умножая выражение (9) слева на T после некоторых несложных преобразований можно доказать, что решение задачи Коши (1) имеет следующее представление

$$x(t)=\sum_{s=0}^{2q}A^{\frac{s-2q}{2p+1}}\frac{d}{dt}\int_{t_0}^{t}\frac{(t-\tau)^{\frac{s-2q}{2q+1}}}{\frac{(s-2q)}{2q+1}!}e^{(\tau-t_0)A^{\frac{2q+1}{2p+1}}}d\tau\,x(t_0). \tag{10}$$

Таким образом, получим следующее утверждение.

**Теорем 1.** Пусть постоянная квадратная матрица A из (1) имеет различные вещественные собственные значения $\lambda_i$. Тогда решение задачи Коши (1) представляется в виде (10).

Отделяя (10) на s=2q и сумму от s=0 до s=2q-1 имеем

$$x(t)=\sum_{s=1}^{2q-1}A^{\frac{s-2q}{2p+1}}\frac{d}{dt}\int_{t_0}^{t}\frac{(t-\tau)^{\frac{s-2q}{2q+1}}}{\frac{(s-2q)}{2q+1}!}e^{(\tau-t_0)A^{\frac{2q+1}{2p+1}}}d\tau\,x(t_0)+\frac{d}{dt}\int_{t_0}^{t}e^{(\tau-t_0)A^{\frac{2q+1}{2p+1}}}d\tau\,x(t_0)=$$

$$=\sum_{s=1}^{2q-1}A^{\frac{s-2q}{2p+1}}\frac{d}{dt}\int_{t_0}^{t}\frac{(t-\tau)^{\frac{s-2q}{2q+1}}}{\frac{s-2q}{2q+1}!}e^{(\tau-t_0)A^{\frac{2q+1}{2p+1}}}d\tau\,x(t_0)+e^{(t-t_0)A^{\frac{2q+1}{2p+1}}}x(t_0) \tag{10$'$}$$

где, при $p=0$, $q=0$ (10$'$) совпадает с известным классическим решением.[7,10]

**Замечание 1.** В случае когда $\alpha=1$, т.е. при $p=0$, $q=0$ из (10) получаем что $s=0$ и $x(t)=\dfrac{d}{dt}\int_{t_0}^{t}e^{A(\tau-t_0)}d\tau x(t_0)=e^{A(t-t_0)}x(t_0)$, которое совпадает с классическим решением задачи Коши (1) при $\alpha=1$.

## 4. Случай кратных или комплексных собственных значений.

Сначала рассмотрим случай, когда среди собственных значений матрица A имеет кратные и комплексные. Заменим A матрицей $A+\varepsilon B$ так, чтобы полученная матрица имела различные вещественные собственные значения $\lambda_i(\varepsilon)$. Тогда решение представляется в виде (10). Учитывая непрерывную зависимость решения от параметра $\varepsilon$, переходя к пределу при $\varepsilon \to 0$, получим результат для кратных или же комплексных собственных значений [11].

**Замечание 2**. Соотношение (10) для решения задачи Коши (1) верно и для случая, когда среди собственных значений матрицы A имеются кратные или комплексные.

## 5. Вычислительный алгоритм.

Более интересным является предложение вычислительной процедуры нахождения $x(t)$ из (8) в скалярном случае, или матричном - из (10). Отметим, что в (9) провести дифференцирование интеграла невозможно из-за отрицательной степени $(t-\tau)$, входящей в подынтегральное выражение. Аналогично, такого рода трудности встречается и в (10). Поэтому, сперва рассмотрим скалярный случай (8) и образуем интеграл из (8) так, чтобы степень $(t-\tau)$ была преобразована в положительную. Для этого отдельно рассмотрим подынтегральное выражение в (8) и его представим в виде:

$$\frac{d}{dt}\int_{t_0}^{t}\frac{(t-\tau)^{\frac{s-2q}{2q+1}}}{\frac{s-2q}{2q+1}!}e^{(\tau-t_0)\lambda_i^{\frac{2q+1}{2p+1}}}\cdot d\tau = -\frac{d}{dt}\int_{t_0}^{t}e^{(\tau-t_0)\lambda_i^{\frac{2q+1}{2p+1}}}d\left(\frac{(t-\tau)^{\frac{s+1}{2q+1}}}{\frac{s+1}{2q+1}!}\right)$$

Далее интегрируя последнее по частям, имеем:

$$-\frac{d}{dt}\left[\frac{(t-\tau)^{\frac{s+1}{2q+1}}}{\frac{s+1}{2q+1}!}e^{(\tau-t_0)\lambda_i^{\frac{2q+1}{2p+1}}}\bigg|_{\tau=t_0}^{t} - \lambda_i^{\frac{2q+1}{2p+1}}\int_{t_0}^{t}\frac{(t-\tau)^{\frac{s+1}{2q+1}}}{\frac{s+1}{2q+1}!}e^{(\tau-t_0)\lambda_i^{\frac{2q+1}{2p+1}}}d\tau\right]=$$

$$=-\frac{d}{dt}\left[-\frac{(t-t_0)^{\frac{s+1}{2q+1}}}{\frac{s+1}{2q+1}!}-\lambda_i^{\frac{2q+1}{2p+1}}\int_{t_0}^{t}\frac{(t-\tau)^{\frac{s+1}{2q+1}}}{\frac{s+1}{2q+1}!}e^{(\tau-t_0)\lambda_i^{\frac{2q+1}{2p+1}}}d\tau\right]= \quad (11)$$

$$=\frac{(t-t_0)^{\frac{s-2q}{2q+1}}}{\frac{s-2q}{2q+1}!}+\lambda_i^{\frac{2q+1}{2p+1}}\int_{t_0}^{t}\frac{(t-\tau)^{\frac{s-2q}{2q+1}}}{\frac{s-2q}{2q+1}!}e^{(\tau-t_0)\lambda_i^{\frac{2q+1}{2p+1}}}d\tau.$$

Подставляя (11) в (8) имеем

$$y_i(t)=\sum_{s=0}^{2q}\lambda_i^{\frac{s-2q}{2p+1}}\left[\frac{(t-t_0)^{\frac{s-2q}{2q+1}}}{\frac{s-2q}{2q+1}!}+\lambda_i^{\frac{2q+1}{2p+1}}\int_{t_0}^{t}\frac{(t-\tau)^{\frac{s-2q}{2q+1}}}{\frac{s-2q}{2q+1}!}e^{(\tau-t_0)\lambda_i^{\frac{2q+1}{2p+1}}}d\tau\right]y_i(t_0) \quad (12)$$

Легко видно из последнего выражения (12), что для y(t) имеется следующая формула

$$y(t) = \sum_{s=0}^{2q} \Lambda^{\frac{s-2q}{2p+1}} \left[ \frac{(t-t_0)^{\frac{s-2q}{2q+1}}}{\frac{s-2q}{2q+1}!} E + \Lambda^{\frac{2q+1}{2p+1}} \int_{t_0}^{t} \frac{(t-\tau)^{\frac{s-2q}{2q+1}}}{\frac{s-2q}{2q+1}!} e^{(\tau-t_0)\Lambda^{\frac{2q+1}{2p+1}}} d\tau \right] y(t_0) \qquad (13)$$

где Е-единичная матрица с соответствующей размерности.

Умножив (13) слева на T и после некоторых преобразований получим следующее выражение для решения задачи Коши (1)

$$x(t) = \sum_{s=1}^{2q-1} A^{\frac{s-2q}{2p+1}} \left[ \frac{(t-t_0)^{\frac{s-2q}{2q+1}}}{\frac{s-2q}{2q+1}!} E + A^{\frac{2q+1}{2p+1}} \int_{t_0}^{t} \frac{(t-\tau)^{\frac{s-2q}{2q+1}}}{\frac{s-2q}{2q+1}!} e^{(\tau-t_0)A^{\frac{2q+1}{2p+1}}} d\tau \right] x(t_0) +$$

$$+ \left[ E + A^{\frac{2q+1}{2p+1}} \int_{t_0}^{t} e^{(\tau-t_0)A^{\frac{2q+1}{2p+1}}} d\tau \right] x(t_o) =$$

$$= \sum_{s=1}^{2q-1} A^{\frac{s-2q}{2p+1}} \left[ \frac{(t-t_0)^{\frac{s-2q}{2q+1}}}{\frac{s-2q}{2q+1}!} E + A^{\frac{2q+1}{2p+1}} \int_{t_0}^{t} \frac{(t-\tau)^{\frac{s-2q}{2q+1}}}{\frac{s-2q}{2q+1}!} e^{(\tau-t_0)A^{\frac{2q+1}{2p+1}}} d\tau \right] x(t_0) + e^{(\tau-t_0)A^{\frac{2q+1}{2p+1}}} x(t_0) \qquad (14)$$

Теперь выражение (12 (или(14)) поддаются использованию для вычислений (в обоих случаях уже отсутствует дифференцирование). Действительно, используя обычные правила прямоугольников [12] для вычислениит $x(t)$ в каждой точке $t_k$ из промежутка $[t_0, T]$ решение (12) (или (14)) имеет вид:

$$y_i(t_k) = \sum_{s=1}^{2q-1} \lambda_i^{\frac{s-2q}{2p+1}} \left[ \frac{(t_k-t_0)^{\frac{s-2q}{2q+1}}}{\frac{s-2q}{2q+1}!} + \lambda_i^{\frac{2q+1}{2p+1}} \sum_{\sigma=0}^{k-1} \frac{(t_k-t_\sigma)^{\frac{s-2q}{2q+1}}}{\frac{s-2q}{2q+1}!} (t_{\sigma+1}-t_\sigma) e^{(t_\sigma-t_0)\lambda_i^{\frac{2q+1}{2p+1}}} \right] y_i(t_0) +$$

$$+ e^{(t_k-t_0)\lambda^{\frac{2q+1}{2p+1}}} y_i(t_0) \qquad (15)$$

Аналогично, для общего случая имеем подобное (15) следующую формулу

$$x(t_k) = \sum_{s=1}^{2q-1} A^{\frac{s-2q}{2p+1}} \left[ \frac{(t_k-t_0)^{\frac{s-2q}{2q+1}}}{\frac{s-2q}{2q+1}!} E + A^{\frac{2q+1}{2p+1}} \sum_{\sigma=0}^{k-1} \frac{(t_k-t_\sigma)^{\frac{s-2q}{2q+1}}}{\frac{s-2q}{2q+1}!} e^{(t_\sigma-t_0)A^{\frac{2q+1}{2p+1}}} (t_{\sigma+1}-t_\sigma) \right] x(t_0) +$$

$$+ e^{(t_k-t_0)A^{\frac{2q+1}{2p+1}}} x(t_0), \qquad (16)$$

Теперь проиллюстрируем выше полученные результаты на следующем примере.

**Пример.** Сначала рассмотрим следующий числовой пример задачи Коши в скалярном случае. Пусть

$$D^\alpha x = a\, x, \quad x(t_0) = x_0 \qquad (17)$$

Применяя формулы (15) для решения задачи (17) численно, покажем, что при $\alpha \to 1$ полученные выражения стремятся к $e^{at}$. В данном примере рассмотрен интервал [0.01  1.01], который разбит шагами на h=0.01 и h=0.001. Такой численный эксперимент показывает правдоподобность полученных формул (8)-(10), (10'). Полученные числовые результаты показывает, что при $\alpha \to 1$ решение уравнения (17) приближается к $e^{2t}$ с погрешностью $nev = \left\| D^\alpha x - \dot{x} \right\|$.

В таблице1 показаны результаты вычисления.

Таблица1. Численная реализация формулы (15)

| a | α | [0.01 1.01] h=0.01 nev | [0.01 1.01] h=0.001 nev |
|---|---|---|---|
| 2 | | | |
| | 1/3 | 1.139e+15 | 3.17e+15 |
| | 3/7 | 4.821e+3 | 1.48e+4 |
| | 199/203 | 1.937e+2 | 6.052e+2 |
| | 1999/2003 | 1.86e+2 | 5.82e+2 |
| | 1 | 2.11e-14 | 7.14e-13 |
| -2 | | | |
| | 1/3 | 19.41 | |
| | 3/7 | 23.59 | |
| | 199/203 | 24.338 | |
| | 1999/2003 | 24.332 | |
| | 1 | 1.11e-15 | |

Теперь иллюстрируем результаты из табл.1 графически на фиг.1 для решения (17) при различных значениях α.

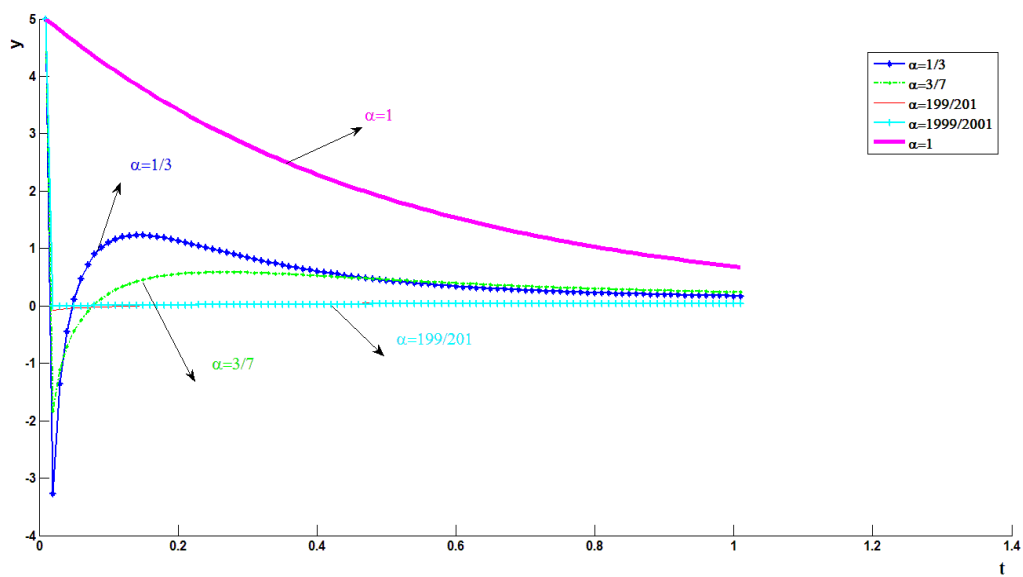

Фиг.1 a=-2 h=0.01, [0.01 1.01]

Из полученных результатов таблицы1 видно, что, погрешность на много больше а это наверное связано с вычислением интеграла методом прямоугольников.

В таблице2 вместо формулы (15), используя формулы (10) и при вычислении интеграла, применяя метод Симпсона т.е. процедуру quad.m пакета МАТЛАБ, которая улучшает полученные результаты.

Таблица2. Численная реализация с методом Симпсона при a=-2

| α | [0.01 1.01] h=0.01 nev | [0.001 10.001] nev |
|---|---|---|
| 1/3 | 12.47 | 4.88 |
| 3/7 | 21.29 | 4.85 |
| 199/203 | 8.4035 | 0.7 |
| 1999/2003 | 2.6539 | 0.5 |
| 1 | 0.6 | 0.002 |

В фиг.2 демонстрируется график функции таблицы2. при различных значениях □

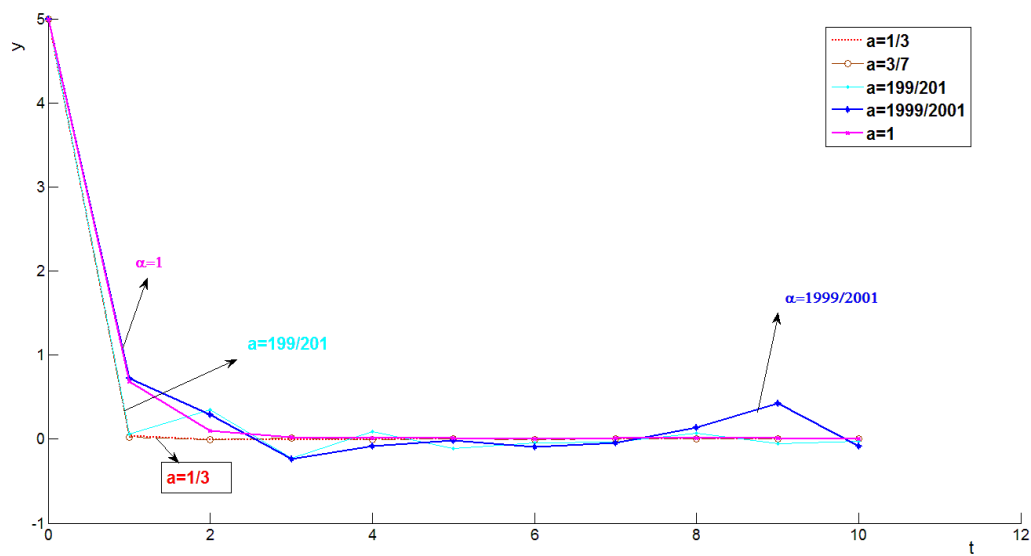

Фиг.2. a=-2, h=1, [1.01 10.01]

Таким образом, при $\alpha = 1$ решение (17) тоже сходится к $e^{-2t}$ и при решении $t \to \infty$ (17) стремится к нулю, т.е., решение последнего уравнения является асимптотически устойчивым.

**Заключение.**

Приводится новая формула для решения задачи Коши системы обыкновенных линейных дробно-производных дифференциальных уравнений с постоянными матричными коэффициентами. Показывается, что эта формула справедлива, когда матричный коэффициент имеет собственные значения как различные, так и кратные или комплексные числа. Далее, полученные формулы с помощью соответствующих преобразований приводятся к виду, который можно использовать для вычислительных целей. Отметим, что тоже представление позволяет обобщить полученные результаты для линейно квадратичных задач оптимального управления в конечном [13-16] и бесконечном интервалах времени [17-22] как в непрерывном, так и в дискретном случаях.